\documentclass[12pt]{article}
\usepackage{latexsym,amsfonts,amsthm,amsmath,amscd,amssymb}
\usepackage[dvips]{graphicx}

\newtheorem{lemma}{Lemma}[section]
\newtheorem{thm}[lemma]{Theorem}
\newtheorem{rem}[lemma]{Remark}
\newtheorem{prop}[lemma]{Proposition}

\newcommand\matS{{\mathbb{S}}}
\newcommand\matP{{\mathbb{P}}}

\newcommand\matZ{{\mathbb{Z}}}

\renewcommand{\hbar}{{\overline{h}}}

\newfont{\Got}{eufm10 scaled 1200}

\newcommand{\mycap} [1] {\caption{\footnotesize{#1}}}

\newcommand{\calO}{{\mathcal{O}}}
\newcommand{\calM}{{\mathcal{M}}}

\newcommand{\hyp}{{\mathrm{hyp}}}
\newcommand{\non}{{\mathrm{non}}}

\begin{document}

\title{The 191 orientable octahedral manifolds}

\author{Damian~\textsc{Heard}\and
Ekaterina~\textsc{Pervova}\and
Carlo~\textsc{Petronio}}

\maketitle

\begin{abstract}
\noindent We enumerate all spaces obtained by gluing in pairs the
faces of the octahedron in an orientation-reversing fashion.
Whenever such a gluing gives rise to non-manifold points, we
remove small open neighbourhoods of these points, so we actually
deal with three-dimensional manifolds with (possibly empty)
boundary.

There are 298 combinatorially inequivalent gluing patterns, and we
show that they define 191 distinct manifolds, of which 132 are
hyperbolic and 59 are not. All the 132 hyperbolic manifolds were
already considered in different contexts by other authors, and we
provide here their known ``names'' together with their main
invariants. We also give the connected sum and JSJ decompositions
for the 59 non-hyperbolic examples.

Our arguments make use of tools coming from hyperbolic geometry,
together with quantum invariants and more classical techniques based on essential
surfaces. Many (but not all) proofs were carried out by computer,
but they do not involve issues of numerical accuracy.

\end{abstract}

\noindent At the very beginning of his fundamental
book~\cite{bible}, as an example of the richness of topology in
three dimensions, Bill Thurston mentioned the fact that there are
quite a few inequivalent ways of gluing together in pairs the
faces of the octahedron. However, to our knowledge, as of today nobody had
ever exactly determined the number of non-homeomorphic 3-manifolds
arising as the results of these gluings. In this note we give
a full solution to this problem, in the context of orientable
(but unoriented) manifolds. After proving that there are 298
inequivalent gluing patterns, we have in fact proved the following:

\begin{thm}\label{main:teo}
Let $O$ be the octahedron and let $\calO$ be the set of homeomorphism
types of $3$-manifolds that can be obtained as follows:
\begin{itemize}
\item First, glue together in pairs in a simplicial and orientation-reversing
fashion the faces of $O$, thus getting a compact polyhedron $X$;
\item Second, remove from $X$ disjoint open stars of the non-manifold points, thus getting
a compact orientable $3$-manifold with (possibly empty) boun\-dary,
all the components of which have positive genus.
\end{itemize}
Then $\calO$ contains precisely $191$ elements, of which $132$ are hyperbolic and $59$ are not.
More precisely, the numbers of inequivalent gluings and manifolds are split according
to the topological type of the boundary as shown in Table~\ref{final:tab}, where $\Sigma_g$
denotes the orientable surface of genus $g$.
\end{thm}

\begin{table}
\begin{center}
\begin{tabular}{c||c|c|c||c}
boundary type & $\#$(gluings) & hyperbolic & non-hyperbolic & total
\\ \hline\hline $\emptyset$ &  37 & -- & 17 & 17\\ \hline $\Sigma_1$ &  81
& $9$ & $21$ & 30\\ \hline $\Sigma_1\sqcup \Sigma_1$ & 9 & 2 & 5 & 7\\ \hline
$\Sigma_2$ &  113 & 63 & 16 & 79\\ \hline $\Sigma_2\sqcup \Sigma_1$ & 2 &
2 & -- & 2 \\ \hline $\Sigma_3$  & 56 & 56  & -- & 56 \\
\hline\hline
Total & 298 &  132 &  59 &  191  \\
\end{tabular}
\end{center}
\mycap{Numbers of distinct manifolds arising from
orientation-reversing gluings of the faces of an
octahedron, with open stars of the non-manifold
points removed\label{final:tab}}
\end{table}

For the PL notions of \emph{polyhedron}, \emph{manifold} and \emph{star}, see for instance~\cite{RS}.
As usual~\cite{lectures,ratcliffe,bible} a 3-manifold $M$ is ``hyperbolic'' if $M$
minus the boundary components of $M$ homeomorphic to the torus carries a complete metric
with constant sectional curvatures $-1$ and totally geodesic boundary.
The removed tori give rise to the so-called \emph{cusps} of the manifold.

\bigskip

In addition to proving Theorem~\ref{main:teo}, we provide rather
detailed information on the 191 elements of $\calO$. In particular,
we determine the volume and other invariants for the 132 hyperbolic
manifolds in $\calO$, and we identify the ``names'' they were given
either in the Callahan-Hildebrand-Weeks census~\cite{CaHiWe,SnapPea}
of small cusped hyperbolic manifolds, or in the
Frigerio-Martelli-Petronio census~\cite{fmp3,www:CP} of small
hyperbolic manifolds with geodesic boundary. We also give detailed
descriptions for the 59 non-hyperbolic elements of $\calO$.

\bigskip

The question of counting the elements of $\calO$ has a rather
transparent combinatorial flavour and appears to be well-suited to
computer investigation, but the complete answer would be extremely
difficult to obtain without the aid of some rather sophisticated
geometric tools developed over the last three decades by a number of
mathematicians. It is indeed mostly thanks to hyperbolic geometry
that one is able to show that certain gluings of $O$, despite being
very similar to each other under many respects, are in fact
distinct. This can be viewed as a manifestation of the crucial
r\^ole played by hyperbolic geometry in the context of
three-dimensional topology, as chiefly witnessed by Thurston's
geometrization, now apparently proved
by Perelman~\cite{bible,perelman1,perelman2,perelman3}

\bigskip

To prove Theorem~\ref{main:teo} we have written some small specific
Haskell code (to list the combinatorially inequivalent gluing
patterns), and then we have used the ``Orb'' and ``Manifold
Recognizer'' programs~\cite{Orb,recognizer}.
There were however some manifolds the computer was unable to
find hyperbolic structures for, and some pairs of manifolds that it was
unable to tell apart.  In these instances, we had to work by hand using classical
techniques, including properly embedded essential surfaces.
Despite being based on computers, our arguments do not involve
issues of numerical accuracy, because approximation was only used
within ``Orb,'' but the results were later checked through exact
arithmetic in algebraic numbers fields with the program
``Snap''~\cite{Snap}.

\vspace{1cm}

\noindent\textsc{Acknowledgements}. Part of this work was carried
out while the third-named author was visiting the University of
Melbourne, the Universit\'e Paul Sabatier in Toulose and the
Columbia University in New York. He is grateful to all these
institutions for financial support, and he would like to thank Craig
Hodgson, Michel Boileau and Dylan Thurston for their warm
hospitality and inspiring mathematical discussions. The second named
author was supported by the Marie Curie fellowship
MIF1-CT-2006-038734.

\section{Preliminaries}

In this section we collect some elementary facts needed to prove
Theorem~\ref{main:teo}.

\paragraph{Polyhedra vs manifolds}
Given a gluing pattern $\varphi$ for the faces of the octahedron $O$, as described in the statement
of Theorem~\ref{main:teo}, let us denote by $X(\varphi)$ the polyhedron
resulting from the gluing, and by $M(\varphi)$ the 3-manifold obtained from $X(\varphi)$
by removing disjoint open stars of the non-manifold points. The following easy fact,
that we leave to the reader, shows that $X(\varphi)$ and $M(\varphi)$ are in fact
very tightly linked:

\begin{prop}
\begin{itemize}
\item Only the points of $X(\varphi)$ arising from the vertices of $O$ can
be non-manifold points of $X(\varphi)$;
\item The homeomorphism type of $X(\varphi)$ determines that of $M(\varphi)$, and
conversely.
\end{itemize}
\end{prop}

Before proceeding, recall that $\calO$ denotes the set of
homeomorphism classes of all $M(\varphi)$'s as
$\varphi$ varies in the set simplicial and orientation-reversing gluing patterns of
the faces of $O$.

\paragraph{Number of inequivalent gluings}
To count the elements of $\calO$, the first step is of course to enumerate the combinatorially
inequivalent gluing patterns $\varphi$. Since $O$ has 8 faces and there are 3 different
ways of gluing together any two chosen faces, the number of different patterns is
$(8-1)!!\times 3^4=105\times 81=8505$. However there is a symmetry group
with 48 elements acting on $O$, so the inequivalent patterns
are actually much fewer than 8505. Using a small piece of Haskell code we have
in fact shown the following:

\begin{prop}\label{patterns:prop}
There exist $298$ combinatorially inequivalent patterns of
orientation-reversing gluings of the faces of $O$.
\end{prop}

\paragraph{Classification according to boundary type}
Two homeomorphic manifolds of course have homeomorphic boundaries.
Moreover the boundary of an orientable 3-manifold is an orientable
surface, which is very easy to identify by counting the number of
connected components and computing the Euler characteristic of each
of them. So the first easy step towards
understanding $\calO$ and proving Theorem~\ref{main:teo} is to split
the inequivalent gluing patterns according to the boundary they give
rise to. Using again a Haskell program we found the results
described in Table~\ref{bd:tab}, where again $\Sigma_g$ denotes the
orientable surface of genus $g$.

\begin{table}
\begin{center}
\begin{tabular}{c|c}
$\partial  M(\varphi)$ & $\#$(inequivalent $\varphi$'s) \\
\hline\hline $\emptyset$ &  37 \\ \hline $\Sigma_1$ & 81\\ \hline $\Sigma_1\sqcup
\Sigma_1$ & 9\\ \hline $\Sigma_2$ &  113\\ \hline $\Sigma_2\sqcup \Sigma_1$ &  2\\
\hline $\Sigma_3$  & 56  \\ \hline\hline
Total &  298 \\
\end{tabular}
\end{center}
\mycap{Numbers of inequivalent gluings $\varphi$ of the faces of
$O$, subdivided according to the topological type of $\partial
 M(\varphi)$\label{bd:tab}}
\end{table}

\paragraph{Further notation}
Choosing one representative for each equivalence class of gluing
patterns $\varphi$ and constructing the corresponding manifold
$M(\varphi)$, we get a set of 298 manifolds that we denote
henceforth by $\calM$. By definition, $\calO$ is obtained from
$\calM$ by identifying homeomorphic manifolds, and the main issue in
establishing Theorem~\ref{main:teo} is indeed to determine which
elements of $\calM$ are in fact homeomorphic to each other. Taking
advantage of the easy work already described, we denote by
$\calM_\Sigma$ the set of elements of $\calM$ having boundary
$\Sigma$, thus getting a splitting of $\calM$ as
$$\calM=
\calM_\emptyset\sqcup
\calM_{\Sigma_1}\sqcup
\calM_{\Sigma_1\sqcup\Sigma_1}\sqcup
\calM_{\Sigma_2}\sqcup
\calM_{\Sigma_2\sqcup\Sigma_1}\sqcup
\calM_{\Sigma_3}.$$
Each set $\calM_\Sigma$, after identifying homeomorphic manifolds,
gives rise to some $\calO_\Sigma$, that we further split as
$$\calO_\Sigma=\calO^\hyp_\Sigma\sqcup\calO^\non_\Sigma,$$
separating the hyperbolic members from the non-hyperbolic ones.

\section{Hyperbolic manifolds}

According to the well-known \emph{rigidity} theorem~\cite{bible,lectures,ratcliffe},
each $3$-manifold carries, up to isometry, at most one
hyperbolic structure, as defined after the statement of Theorem~\ref{main:teo}.
Note that the hyperbolic structures we consider are finite-volume by default.
Moreover the following facts hold:

\begin{enumerate}

\item Every hyperbolic manifold
with cusps or non-empty boundary
has a ``canonical decomposition,'' which allows
to efficiently compare it to any other such manifold. This is the decomposition
into ideal polyhedra due to Epstein and Penner~\cite{EpPe} for \emph{cusped} manifolds
(non-compact and without boundary), and the decomposition into truncated hyperideal polyhedra
due to Kojima~\cite{koji1,koji2} for manifolds with non-empty boundary.
The hyperbolic structure of a manifold, whence (by rigidity) its topology, determines
not only the polyhedral type of the blocks of the decomposition, but
also the combinatorics of the gluings;

\item If a manifold is
represented by a \emph{triangulation}, namely as a gluing of tetrahedra,
both its hyperbolic structure (if any) and its canonical decomposition can be searched algorithmically.
This applies in particular to any element of the set $\calM$ of manifolds we need to
analyze, because the octahedron $O$ can be viewed as a partial gluing of 4 tetrahedra.
The idea to construct the hyperbolic structure, due to Thurston~\cite{bible},
is to consider a space of parameters for the hyperbolic structures on each individual
tetrahedron, and then to express the matching of the structures on the glued tetrahedra
by a system of equations, that can be solved using numerical tools.
The method for constructing the canonical decomposition is to modify any given geometric
triangulation until the canonical decomposition is reached. This uses
the ``tilt formula'' of Sakuma and Weeks~\cite{SaWe} for cusped
manifolds, and its variation due to Ushijima~\cite{akira:deco}, together with
some ideas from~\cite{fp}, for manifolds with non-empty geodesic boundary.
Both the search for the hyperbolic structure and that for the canonical decomposition
are not fully guaranteed to work, but in practice they always do (perhaps after
some initial randomization of the triangulation);

\item The computer programs ``SnapPea''~\cite{SnapPea} by Jeff Weeks, and
``Orb''~\cite{Orb} by Damian Heard very efficiently implement the procedures
mentioned in the previous point;

\item Both ``SnapPea'' and ``Orb'' employ numerical approximation, but the solutions
these programs find can be checked using exact arithmetic in
algebraic number fields with the
program ``Snap''~\cite{Snap} by Oliver Goodman.

\end{enumerate}

\paragraph{Genus-3 geodesic boundary}
To prove Theorem~\ref{main:teo}, for each of the 6 sets $\calM_\Sigma$ we have,
we need to determine which elements of $\calM_\Sigma$ are homeomorphic to each other,
thus finding the corresponding $\calO_\Sigma$, and then to decide which elements of
$\calO_\Sigma$ are actually hyperbolic. We begin with the case $\Sigma=\Sigma_3$,
where the result is quite striking. It
was initially discovered from a computer experiment~\cite{fmp3} and
later established theoretically. We include a sketch of
the proof for the sake of completeness.

\begin{prop}\label{56:prop}
The $56$ elements of $\calM_{\Sigma_3}$ are all hyperbolic and distinct from each other,
so $\calO_{\Sigma_3}=\calO_{\Sigma_3}^\hyp=\calM_{\Sigma_3}$ has $56$ elements. For each
element of this set the
Kojima canonical decomposition has the same single block, namely a truncated regular
hyperbolic octahedron with all dihedral angles equal to $\pi/6$.
\end{prop}

\begin{proof}
An easy computation of Euler characteristic shows that a gluing
$\varphi$ defines a manifold $ M(\varphi)$ bounded by $\Sigma_3$
if and only if it identifies all 12 edges to each other. We want to
show that such an $ M(\varphi)$ is hyperbolic with geodesic
boundary by choosing a hyperbolic shape of the truncated octahedron
that is matched by $\varphi$. Since all edges are glued
together, this can only happen if the geometric shape is such that
all edges have the same length, \emph{i.e.}~the octahedron is regular.
If this is the case, all dihedral angles are also the same, so they
must all be $2\pi/12$. Such an octahedron certainly does not exist
in Euclidean or spherical geometry, but it does in hyperbolic
geometry. This implies that $ M(\varphi)$ is indeed hyperbolic.

Let us now analyze the Kojima canonical decomposition of
$ M(\varphi)$. To this end we recall~\cite{koji1,koji2} that it
is dual to the cut locus of the boundary, \emph{i.e.}~to the set of
points having multiple shortest paths to $\partial  M(\varphi)$.
Using the fact that $ M(\varphi)$ is the gluing of a regular
truncated octahedron, which is totally symmetric, it is not too
difficult to show that the Kojima decomposition is given by the
octahedron itself, with its gluing pattern $\varphi$. This implies
that the geometry of $ M(\varphi)$, and hence its topology,
determines $\varphi$. Therefore different $\varphi$'s give rise to
different $ M(\varphi)$'s.
\end{proof}

It follows from this result that the 56 elements of
$\calO_{\Sigma_3}^\hyp$ all have the same volume, that one can
calculate to be $11.448776110...$ via Ushijima's
formulae~\cite{akira:vol}. Using ``Orb'' we have also computed the
symmetry groups and homology of the elements of
$\calO_{\Sigma_3}^\hyp$, as described in Tables~\ref{bdS3:tab:1}
and~\ref{bdS3:tab:2}. Note that these invariants alone are far from
sufficient to distinguish the 56 elements of $\calM_{\Sigma_3}$ from
each other. The tables also show the position of the manifolds in
the file {\tt census$\_$4$\_$T3$\_$octa.snp} available
from~\cite{www:CP}. Here and below $\matZ_n$ and $D_n$ denote
respectively the cyclic group with $n$ elements and the dihedral
group with $2n$ elements.

\begin{table}
\begin{center}
\begin{tabular}{l|c|c|c|c}
File & no. & Volume & Sym & Hom \\
\hline \hline {\tt census$\_$4$\_$T3$\_$octa.snp} &    0  &
11.448776110 & trivial & $\matZ^3$ \\ \hline {\tt
census$\_$4$\_$T3$\_$octa.snp} &    1  &  11.448776110 & trivial &
$\matZ^3$ \\ \hline {\tt census$\_$4$\_$T3$\_$octa.snp} &    2  &
11.448776110 & trivial & $\matZ^3$ \\ \hline {\tt
census$\_$4$\_$T3$\_$octa.snp} &    3  &  11.448776110 & trivial &
$\matZ^3$ \\ \hline {\tt census$\_$4$\_$T3$\_$octa.snp} &    4  &
11.448776110 & trivial & $\matZ^3$ \\ \hline {\tt
census$\_$4$\_$T3$\_$octa.snp} &    5  &  11.448776110 & trivial &
$\matZ^3$ \\ \hline {\tt census$\_$4$\_$T3$\_$octa.snp} &    6  &
11.448776110 & trivial & $\matZ^3$ \\ \hline {\tt
census$\_$4$\_$T3$\_$octa.snp} &    7  &  11.448776110 & trivial &
$\matZ^3$ \\ \hline {\tt census$\_$4$\_$T3$\_$octa.snp} &    8  &
11.448776110 & trivial & $\matZ^3$ \\ \hline {\tt
census$\_$4$\_$T3$\_$octa.snp} &    9  &  11.448776110 & trivial &
$\matZ^3$ \\ \hline {\tt census$\_$4$\_$T3$\_$octa.snp} &   10  &
11.448776110 & trivial & $\matZ^3$ \\ \hline {\tt
census$\_$4$\_$T3$\_$octa.snp} &   11  &  11.448776110 & trivial &
$\matZ^3$ \\ \hline {\tt census$\_$4$\_$T3$\_$octa.snp} &   12  &
11.448776110 & $D_2$ & $\matZ^3$ \\ \hline {\tt
census$\_$4$\_$T3$\_$octa.snp} &   13  &  11.448776110 & $\matZ_2$ &
$\matZ^3$ \\ \hline {\tt census$\_$4$\_$T3$\_$octa.snp} &   14  &
11.448776110 & $\matZ_2$ & $\matZ^3$ \\ \hline {\tt
census$\_$4$\_$T3$\_$octa.snp} &   15  &  11.448776110 & trivial &
$\matZ^3$ \\ \hline {\tt census$\_$4$\_$T3$\_$octa.snp} &   16  &
11.448776110 & $D_2$ & $\matZ^3$ \\ \hline {\tt
census$\_$4$\_$T3$\_$octa.snp} &   17  &  11.448776110 & $\matZ_2$ &
$\matZ^3$ \\ \hline {\tt census$\_$4$\_$T3$\_$octa.snp} &   18  &
11.448776110 & $D_4$ & $\matZ^3$ \\ \hline {\tt
census$\_$4$\_$T3$\_$octa.snp} &   19  &  11.448776110 & $\matZ_2$ &
$\matZ^3$ \\ \hline {\tt census$\_$4$\_$T3$\_$octa.snp} &   20  &
11.448776110 & $\matZ_2$ & $\matZ^3$ \\ \hline {\tt
census$\_$4$\_$T3$\_$octa.snp} &   21  &  11.448776110 & $D_4$ &
$\matZ^3$ \\ \hline {\tt census$\_$4$\_$T3$\_$octa.snp} &   22  &
11.448776110 & $\matZ_2$ & $\matZ^3$ \\ \hline {\tt
census$\_$4$\_$T3$\_$octa.snp} &   23  &  11.448776110 & $\matZ_2$ &
$\matZ^3$ \\ \hline {\tt census$\_$4$\_$T3$\_$octa.snp} &   24  &
11.448776110 & $\matZ_2$ & $\matZ^3$ \\ \hline {\tt
census$\_$4$\_$T3$\_$octa.snp} &   25  &  11.448776110 & trivial &
$\matZ^3$ \\ \hline {\tt census$\_$4$\_$T3$\_$octa.snp} &   26  &
11.448776110 & trivial & $\matZ^3$ \\ \hline {\tt
census$\_$4$\_$T3$\_$octa.snp} &   27  &  11.448776110 & trivial &
$\matZ^3$ \\ \hline {\tt census$\_$4$\_$T3$\_$octa.snp} &   28  &
11.448776110 & trivial & $\matZ^3$ \\ \hline {\tt
census$\_$4$\_$T3$\_$octa.snp} &   29  &  11.448776110 & trivial &
$\matZ^3$ \\ \hline {\tt census$\_$4$\_$T3$\_$octa.snp} &   30  &
11.448776110 & trivial & $\matZ^3$ \\
\end{tabular}
\end{center}
\mycap{Information on the $56$ elements of $\calO^\hyp_{\Sigma_3}$
(the compact orientable hyperbolic manifolds with geodesic boundary of genus $3$
arising from gluings of the octahedron) -- part 1
\label{bdS3:tab:1}}
\end{table}

\begin{table}
\begin{center}
\begin{tabular}{l|c|c|c|c}
File & no. & Volume & Sym & Hom \\
\hline \hline  {\tt
census$\_$4$\_$T3$\_$octa.snp} &   31  &  11.448776110 & $\matZ_2$ &
$\matZ^3$ \\ \hline {\tt census$\_$4$\_$T3$\_$octa.snp} &   32  &
11.448776110 & $\matZ_2$ & $\matZ^3$ \\ \hline {\tt
census$\_$4$\_$T3$\_$octa.snp} &   33  &  11.448776110 & $\matZ_2$ &
$\matZ^3$ \\ \hline {\tt census$\_$4$\_$T3$\_$octa.snp} &   34  &
11.448776110 & $\matZ_2$ & $\matZ^3$ \\ \hline {\tt
census$\_$4$\_$T3$\_$octa.snp} &   35  &  11.448776110 & $\matZ_2$ &
$\matZ^3$ \\ \hline {\tt census$\_$4$\_$T3$\_$octa.snp} &   36  &
11.448776110 & $\matZ_2$ & $\matZ^3$ \\ \hline {\tt
census$\_$4$\_$T3$\_$octa.snp} &   37  &  11.448776110 & trivial &
$\matZ_3+\matZ^3$ \\ \hline {\tt census$\_$4$\_$T3$\_$octa.snp} & 38
& 11.448776110 & $\matZ_2$ &  $\matZ_3+\matZ^3$ \\ \hline {\tt
census$\_$4$\_$T3$\_$octa.snp} &   39  &  11.448776110 & $D_2$ &
$\matZ_3+\matZ^3$ \\ \hline {\tt census$\_$4$\_$T3$\_$octa.snp} & 40
& 11.448776110 & $\matZ_4$ & $\matZ_3+\matZ^3$ \\ \hline {\tt
census$\_$4$\_$T3$\_$octa.snp} &   41  &  11.448776110 & $\matZ_2$ &
$\matZ^3$ \\ \hline {\tt census$\_$4$\_$T3$\_$octa.snp} &   42  &
11.448776110 & trivial & $\matZ^3$ \\ \hline {\tt
census$\_$4$\_$T3$\_$octa.snp} &   43  &  11.448776110 & $\matZ_2$ &
$\matZ^3$ \\ \hline {\tt census$\_$4$\_$T3$\_$octa.snp} &   44  &
11.448776110 & $\matZ_2$ & $\matZ^3$ \\ \hline {\tt
census$\_$4$\_$T3$\_$octa.snp} &   45  &  11.448776110 & $\matZ_2$ &
$\matZ^3$ \\ \hline {\tt census$\_$4$\_$T3$\_$octa.snp} &   46  &
11.448776110 & trivial & $\matZ^3$ \\ \hline {\tt
census$\_$4$\_$T3$\_$octa.snp} &   47  &  11.448776110 & trivial &
$\matZ^3$ \\ \hline {\tt census$\_$4$\_$T3$\_$octa.snp} &   48  &
11.448776110 & $\matZ_2$ & $\matZ^3$ \\ \hline {\tt
census$\_$4$\_$T3$\_$octa.snp} &   49  &  11.448776110 & trivial &
$\matZ^3$ \\ \hline {\tt census$\_$4$\_$T3$\_$octa.snp} &   50  &
11.448776110 & $\matZ_2$ & $\matZ^3$ \\ \hline {\tt
census$\_$4$\_$T3$\_$octa.snp} &   51  &  11.448776110 & trivial &
$\matZ^3$ \\ \hline {\tt census$\_$4$\_$T3$\_$octa.snp} &   52  &
11.448776110 & trivial & $\matZ^3$ \\ \hline {\tt
census$\_$4$\_$T3$\_$octa.snp} &   53  &  11.448776110 &trivial &
$\matZ^3$ \\ \hline {\tt census$\_$4$\_$T3$\_$octa.snp} &   54  &
11.448776110 & $\matZ_2$ & $\matZ^3$ \\ \hline {\tt
census$\_$4$\_$T3$\_$octa.snp} &   55  &  11.448776110 & $\matZ_2$
& $\matZ^3$ \\
\end{tabular}
\end{center}
\mycap{Information on the $56$ elements of $\calO^\hyp_{\Sigma_3}$
(the compact orientable hyperbolic manifolds with geodesic boundary of genus $3$
arising from gluings of the octahedron) -- part 2
\label{bdS3:tab:2}}
\end{table}

\paragraph{Genus-2 geodesic boundary and one cusp}
For the case $\Sigma=\Sigma_2\sqcup\Sigma_1$ the analysis of $\calM_\Sigma$
was already contained in~\cite{fmp3}:

\begin{prop}\label{Sigma21:prop}
The two elements of $\calM_{\Sigma_2\sqcup\Sigma_1}$ are hyperbolic and
distinct from each other, so
$\calO_{\Sigma_2\sqcup\Sigma_1}=\calO_{\Sigma_2\sqcup\Sigma_1}^\hyp=\calM_{\Sigma_2\sqcup\Sigma_1}$
has two elements.
\end{prop}

Table~\ref{bdS2T:tab} describes the
symmetry group and homology of both elements of $\calO_{\Sigma_2\sqcup\Sigma_1}^\hyp$,
and reference to their position in the files available from~\cite{www:CP}, as we determined using
``Orb''.

\begin{table}
\begin{center}
\begin{tabular}{l|c|c|c|c}
File & no. & Volume & Sym & Hom \\
\hline \hline {\tt census$\_$4$\_$cusp.snp} &  14 &    8.681737155 &
$\matZ_2$ & $\matZ^3$ \\ \hline {\tt census$\_$4$\_$cusp.snp} &  15
& 8.681737155 & $\matZ_2$ & $\matZ^3$ \\
\end{tabular}
\end{center}
\mycap{Information on the $2$ elements of $\calO^\hyp_{\Sigma_2\sqcup\Sigma_1}$
(the orientable hyperbolic manifolds with one cusp and geodesic boundary of genus 2
arising from gluings of the octahedron)
\label{bdS2T:tab}}
\end{table}

\paragraph{Genus-2 geodesic boundary}
The following partial information on the elements of $\calM_{\Sigma_2}$
can be deduced from the results in~\cite{fmp3}:

\begin{prop}\label{FMP:Sigma2:prop}
The set $\calM_{\Sigma_2}$ (which has $113$ elements) contains the following subsets:

\begin{itemize}
\item A set of $14$ distinct hyperbolic manifolds with Kojima decomposition having one
and the same block, namely a regular truncated octahedron with all dihedral angles equal to $\pi/3$;

\item A set of $8$ distinct hyperbolic manifolds with Kojima decomposition having one
and the same block, namely a non-regular truncated octahedron;

\item A set of $4$ distinct hyperbolic manifolds with Kojima decomposition
having the same two blocks, namely two identical square pyramids.

\end{itemize}
Moreover any other hyperbolic element of $\calM_{\Sigma_2}$
has Kojima decomposition consisting of tetrahedra only.
\end{prop}

To complete the analysis of the hyperbolic elements of
$\calM_{\Sigma_2}$, using ``Orb'' (and then ``Snap'' for
a formal verification) we proved the following:

\begin{prop}\label{Damian:Sigma2:prop}
Of the $113-(14+8+4)=87$ elements of $\calM_{\Sigma_2}$
not covered by Proposition~\ref{FMP:Sigma2:prop}, at least $37$ are
hyperbolic, and they are all distinct from each other.
\end{prop}

After ``Orb'' has been able to construct the hyperbolic structure of an element $M$
of $\calM$ and the solution has been checked using ``Snap,'' one can state for sure that
$M$ is indeed hyperbolic, and one can positively determine whether $M$ is homeomorphic
to any other given hyperbolic manifold. However if ``Orb'' fails to construct the structure
one has to prove by some other method that $M$ is actually non-hyperbolic. This is what we
do in the next section. In particular, we prove that
the $113-[(14+8+4)+37]=50$ elements of $\calM_{\Sigma_2}$
not covered by Propositions~\ref{FMP:Sigma2:prop} and~\ref{Damian:Sigma2:prop} are indeed
non-hyperbolic, which implies the following:

\begin{prop}\label{hyp:Sigma2:prop}
The set $\calO^\hyp_{\Sigma_2}$ consists of the $63$ manifolds described
in Propositions~\ref{FMP:Sigma2:prop} and~\ref{Damian:Sigma2:prop}.
\end{prop}

The elements of $\calO^\hyp_{\Sigma_2}$, together with the usual
information on them determined by ``Orb,'' are listed in order of
increasing volume in Tables~\ref{bdS2:tab:1} and~\ref{bdS2:tab:2}.
Again the first column indicates the file from~\cite{www:CP} where
the manifold can be located in the position (starting from 0)
specified in the second column. Note that the name of the file
contains a description of the Kojima canonical decomposition
(\emph{e.g.} {\tt tetra6} means that this decomposition consists of
6 tetrahedra).

\begin{table}
\begin{center}
\begin{tabular}{l|c|c|c|c}
File & no. & Volume & Sym & Hom \\
\hline \hline {\tt census$\_$3.snp} &  93 &    7.636519630 &
$\matZ_2$ & $\matZ^2$ \\ \hline {\tt census$\_$3.snp} &  90 &
7.636519630 & $\matZ_2$ & $\matZ^2$ \\ \hline {\tt census$\_$3.snp}
&  89 &    7.636519630 & $\matZ_2$ & $\matZ^2$ \\ \hline {\tt
census$\_$3.snp} &  88 &    7.636519630 & $\matZ_2$ & $\matZ^2$ \\
\hline {\tt census$\_$3.snp} &  92 &    7.636519630 & trivial &
$\matZ^2$ \\ \hline {\tt census$\_$3.snp} &  86 &    7.636519630 &
trivial & $\matZ^2$ \\ \hline {\tt census$\_$3.snp} &  87 &
7.636519630 & $\matZ_2$ & $\matZ^2$ \\ \hline {\tt census$\_$3.snp}
&  94 &    7.636519630 & $\matZ_2$ & $\matZ^2$ \\ \hline {\tt
census$\_$3.snp} &  91 &    7.636519630 & $\matZ_2$ & $\matZ^2$ \\
\hline {\tt census$\_$4$\_$T2$\_$tetra6.snp} &   2 &    8.297977385
& $\matZ_2$ & $\matZ^2$ \\ \hline {\tt
census$\_$4$\_$T2$\_$tetra6.snp} & 1 &    8.297977385 & $\matZ_2$ &
$\matZ^2$ \\ \hline {\tt census$\_$4$\_$T2$\_$tetra6.snp} &   0 &
8.297977385 & $\matZ_2$ & $\matZ^2$ \\ \hline {\tt
census$\_$4$\_$T2$\_$tetra4.snp} &  75 &    8.625848296 & $\matZ_2$
& $\matZ^2$ \\ \hline{\tt census$\_$4$\_$T2$\_$tetra4.snp} &  76 &
8.625848296 & $\matZ_2$ & $\matZ^2$ \\ \hline {\tt
census$\_$4$\_$T2$\_$octa$\_$nonreg.snp}  &  1 &    8.739252140 &
$D_3$ & $\matZ_3+\matZ^2$ \\ \hline {\tt
census$\_$4$\_$T2$\_$octa$\_$nonreg.snp}  &  0 &    8.739252140 &
$D_3$ & $\matZ_3+\matZ^2$ \\ \hline {\tt
census$\_$4$\_$T2$\_$octa$\_$nonreg.snp}  &  7 &    8.739252140 &
$D_3$ & $\matZ_3+\matZ^2$ \\ \hline {\tt
census$\_$4$\_$T2$\_$octa$\_$nonreg.snp}  &  6 &    8.739252140 &
$D_3$ & $\matZ_3+\matZ^2$ \\ \hline {\tt
census$\_$4$\_$T2$\_$octa$\_$nonreg.snp}  &  2 &    8.739252140 &
$\matZ_2$ & $\matZ^2$ \\ \hline {\tt
census$\_$4$\_$T2$\_$octa$\_$nonreg.snp}  &  5 &    8.739252140 &
$\matZ_2$ & $\matZ^2$ \\ \hline {\tt
census$\_$4$\_$T2$\_$octa$\_$nonreg.snp}  &  4 &    8.739252140 &
$\matZ_2$ & $\matZ^2$ \\ \hline {\tt
census$\_$4$\_$T2$\_$octa$\_$nonreg.snp}  &  3 &    8.739252140 &
$\matZ_2$ & $\matZ^2$ \\ \hline {\tt
census$\_$4$\_$T2$\_$pyramids.snp} & 2 &    9.044841574 & $\matZ_2$
& $\matZ^2$ \\ \hline {\tt census$\_$4$\_$T2$\_$pyramids.snp} &   1
& 9.044841574 & $\matZ_2$ & $\matZ^2$ \\ \hline {\tt
census$\_$4$\_$T2$\_$pyramids.snp} &   0 &    9.044841574 &
$\matZ_2$ & $\matZ^2$ \\ \hline {\tt
census$\_$4$\_$T2$\_$pyramids.snp} & 3 &    9.044841574 & $\matZ_2$
& $\matZ^2$ \\ \hline {\tt census$\_$4$\_$T2$\_$tetra4.snp} & 161  &
9.082538547 & trivial & $\matZ^2$ \\ \hline {\tt
census$\_$4$\_$T2$\_$tetra4.snp} & 162  &   9.082538547 & trivial &
$\matZ^2$ \\ \hline {\tt census$\_$4$\_$T2$\_$tetra4.snp} & 166  &
9.087925790 & $\matZ_2$ & $\matZ^2$ \\ \hline {\tt
census$\_$4$\_$T2$\_$tetra4.snp} & 165  &   9.087925790 & $\matZ_2$
& $\matZ^2$ \\ \hline {\tt census$\_$4$\_$T2$\_$tetra4.snp} & 163  &
9.087925790 & $\matZ_2$ & $\matZ_3+\matZ^2$ \\ \hline {\tt
census$\_$4$\_$T2$\_$tetra4.snp} & 164  &   9.087925790 & $\matZ_2$
& $\matZ_3+\matZ^2$ \\
\end{tabular}
\end{center}
\mycap{Information on the $63$ elements of $\calO^\hyp_{\Sigma_2}$
(the compact orientable hyperbolic manifolds with geodesic boundary of genus $2$
arising from gluings of the octahedron) -- part 1
\label{bdS2:tab:1}}
\end{table}

\begin{table}
\begin{center}
\begin{tabular}{l|c|c|c|c}
File & no. & Volume & Sym & Hom \\
\hline \hline {\tt census$\_$4$\_$T2$\_$tetra5.snp}
& 4 &   9.134474458 & $D_4$ & $\matZ^3$ \\ \hline {\tt
census$\_$4$\_$T2$\_$tetra5.snp} &   3  &   9.134474458 & $\matZ_2$
& $\matZ_2+\matZ^2$ \\ \hline {\tt census$\_$4$\_$T2$\_$tetra5.snp}
& 7 &   9.134474458 & $\matZ_2$ & $\matZ_2+\matZ^2$ \\ \hline {\tt
census$\_$4$\_$T2$\_$tetra5.snp} &   5  &   9.134474458 & $D_4$ &
$\matZ^3$ \\ \hline {\tt census$\_$4$\_$T2$\_$tetra5.snp} &   6  &
9.134474458 & $D_2$ & $\matZ_2+\matZ^2$ \\ \hline {\tt
census$\_$4$\_$T2$\_$tetra5.snp} &   8  &   9.134474458 & $D_2$ &
$\matZ_2+\matZ^2$ \\ \hline {\tt census$\_$4$\_$T2$\_$tetra5.snp} &
15  & 9.333442928 & $\matZ_2$ & $\matZ^2$ \\ \hline {\tt
census$\_$4$\_$T2$\_$tetra5.snp} &  18  &   9.333442928 & $\matZ_2$
& $\matZ^2$ \\ \hline {\tt census$\_$4$\_$T2$\_$tetra5.snp} &  16  &
9.333442928 & trivial & $\matZ^2$ \\ \hline {\tt
census$\_$4$\_$T2$\_$tetra5.snp} &  19  &   9.333442928 & $\matZ_2$
& $\matZ_3+\matZ^2$ \\ \hline {\tt census$\_$4$\_$T2$\_$tetra5.snp}
& 17 &   9.333442928 & trivial & $\matZ^2$ \\ \hline {\tt
census$\_$4$\_$T2$\_$tetra5.snp} &  20  &   9.333442928 & $\matZ_2$
& $\matZ_3+\matZ^2$ \\ \hline {\tt census$\_$4$\_$T2$\_$tetra4.snp}
& 246 &   9.346204962 & trivial & $\matZ^2$ \\ \hline {\tt
census$\_$4$\_$T2$\_$tetra4.snp} & 245  &   9.346204962 & $\matZ_2$
& $\matZ_3+\matZ^2$ \\ \hline {\tt census$\_$4$\_$T2$\_$tetra4.snp}
& 247 &   9.346204962 & $\matZ_2$ & $\matZ^2$ \\ \hline {\tt
census$\_$4$\_$T2$\_$tetra5.snp} &  21  &   9.350261353 & $\matZ_2$
& $\matZ_3+\matZ^2$ \\ \hline {\tt census$\_$4$\_$T2$\_$tetra5.snp}
& 22 &   9.350261353 & $\matZ_2$ & $\matZ^2$ \\ \hline {\tt
census$\_$4$\_$T2$\_$octa$\_$reg.snp} &  11  &   9.415841683 & $D_2$
& $\matZ^2$ \\ \hline {\tt census$\_$4$\_$T2$\_$octa$\_$reg.snp} & 5
& 9.415841683 & $\matZ_2$ & $\matZ/2+\matZ^2$ \\ \hline {\tt
census$\_$4$\_$T2$\_$octa$\_$reg.snp} &   1  &   9.415841683 & $D_4$
& $\matZ/3+\matZ^2$ \\ \hline {\tt
census$\_$4$\_$T2$\_$octa$\_$reg.snp} & 9  &   9.415841683 &
$\matZ_2+\matZ_4$ & $\matZ_5+\matZ^2$ \\ \hline {\tt
census$\_$4$\_$T2$\_$octa$\_$reg.snp} &   6  &   9.415841683 & $D_3$
& $\matZ_6+\matZ^2$ \\ \hline {\tt
census$\_$4$\_$T2$\_$octa$\_$reg.snp} & 7  &   9.415841683 & $D_3$ &
$\matZ_6+\matZ^2$ \\ \hline {\tt
census$\_$4$\_$T2$\_$octa$\_$reg.snp} & 4  &   9.415841683 &
$\matZ_2$ & $\matZ_2+\matZ^2$ \\ \hline {\tt
census$\_$4$\_$T2$\_$octa$\_$reg.snp} &   3  &   9.415841683 &
trivial & $\matZ_2+\matZ^2$ \\ \hline {\tt
census$\_$4$\_$T2$\_$octa$\_$reg.snp} &  10  &   9.415841683 & $D_4$
& $\matZ^2$ \\ \hline {\tt census$\_$4$\_$T2$\_$octa$\_$reg.snp} & 8
& 9.415841683 & trivial & $\matZ_2+\matZ^2$ \\ \hline {\tt
census$\_$4$\_$T2$\_$octa$\_$reg.snp} &  13  &   9.415841683 &
trivial & $\matZ^2$ \\ \hline {\tt
census$\_$4$\_$T2$\_$octa$\_$reg.snp} & 2  &   9.415841683 & $D_2$ &
$\matZ^2$ \\ \hline {\tt census$\_$4$\_$T2$\_$octa$\_$reg.snp} &  12
& 9.415841683 & $\matZ_4$ & $\matZ^2$ \\ \hline {\tt
census$\_$4$\_$T2$\_$octa$\_$reg.snp} &   0  &   9.415841683 &
$D_2$ & $\matZ_3+\matZ^2$ \\
\end{tabular}
\end{center}
\mycap{Information on the $63$ elements of $\calO^\hyp_{\Sigma_2}$
(the compact orientable hyperbolic manifolds with geodesic boundary of genus $2$
arising from gluings of the octahedron) -- part 2
\label{bdS2:tab:2}}
\end{table}

\paragraph{Cusped manifolds}
We carried out the analysis of the hyperbolic elements of $\calM_{\Sigma_1}$
and $\calM_{\Sigma_1\sqcup\Sigma_1}$ using ``Orb,'' with the following result:

\begin{prop}\label{Damian:cusped:prop}
\begin{itemize}
\item The set $\calM_{\Sigma_1}$ (which has $81$ elements)
contains $11$ hyperbolic manifolds,
yielding $9$ distinct homeomorphism types;

\item The set $\calM_{\Sigma_1\sqcup\Sigma_1}$ (which has $9$ elements)
contains $2$ distinct hyperbolic manifolds.
\end{itemize}
\end{prop}

As above for the case of boundary $\Sigma_2$, failure of ``Orb'' to find a cusped hyperbolic
structure does not imply that the structure does not exist. However in the next section
we show that the $81-11=70$ elements of $\calM_{\Sigma_1}$ and the
$9-2=7$ elements of $\calM_{\Sigma_1\sqcup\Sigma_1}$ not covered by
Proposition~\ref{Damian:cusped:prop} are indeed non-hyperbolic, which implies the following:

\begin{prop}\label{hyp:Sigma1:prop}
The set $\calO^\hyp_{\Sigma_1}$ (respectively, $\calO^\hyp_{\Sigma_1\sqcup\Sigma_1}$) consists of the $9$ (respectively, $2$) manifolds described
in Proposition~\ref{Damian:cusped:prop}.
\end{prop}

Using ``Orb'' we have determined the symmetry group and homology of each element of
$\calO^\hyp_{\Sigma_1}$ and $\calO^\hyp_{\Sigma_1\sqcup\Sigma_1}$, together with the name it
was given in~\cite{CaHiWe,SnapPea}. This information appears in
Tables~\ref{bdT:tab} and~\ref{bdTT:tab}.

\begin{table}
\begin{center}
\begin{tabular}{l|c|c|c}
Name & Volume & Sym & Hom \\
\hline \hline m006 & 2.568970601 & $D_2$ & $\matZ_5+\matZ$ \\ \hline
m007 & 2.568970601 & $D_2$ & $\matZ_3+\matZ$ \\ \hline m009 &
2.666744783 & $D_2$ & $\matZ_2+\matZ$ \\ \hline m010 & 2.666744783 &
$D_2$ & $\matZ_6+\matZ$ \\ \hline m011 & 2.781833912 & $\matZ$ &
$\matZ_2$ \\ \hline m032 & 3.163963229 & $D_2$ & $\matZ$ \\ \hline
m033 & 3.163963229 & $D_2$ & $\matZ_9+\matZ$ \\ \hline m036 &
3.177293279 & $D_2$ & $\matZ_3+\matZ$ \\ \hline m038 & 3.177293279 &
$D_2$ & $\matZ$ \\
\end{tabular}
\end{center}
\mycap{Information on the $9$ elements of $\calO^\hyp_{\Sigma_1}$
(the one-cusped orientable hyperbolic manifolds arising from gluings of the octahedron)
\label{bdT:tab}}
\end{table}

\begin{table}
\begin{center}
\begin{tabular}{l|c|c|c}
Name & Volume & Sym & Hom \\
\hline \hline m125 & 3.663862377 & $D_4$ & $\matZ^2$ \\ \hline
m129 & 3.663862377 & $D_4$ & $\matZ^2$ \\
\end{tabular}
\end{center}
\mycap{Information on the $2$ elements of $\calO^\hyp_{\Sigma_1\sqcup\Sigma_1}$
(the two-cusped orientable hyperbolic manifolds arising from gluings of the octahedron)
\label{bdTT:tab}}
\end{table}

\section{Non-hyperbolic manifolds}
In this section we analyze the elements of $\calM$ not covered by
Propositions~\ref{56:prop}, \ref{Sigma21:prop},~\ref{FMP:Sigma2:prop},~\ref{Damian:Sigma2:prop}, and~\ref{Damian:cusped:prop},
thus completing our enumeration of $\calO$. Recall that only $\calM_\emptyset$, $\calM_{\Sigma_1}$,
$\calM_{\Sigma_1\sqcup\Sigma_1}$, and $\calM_{\Sigma_2}$ still require some work.

\paragraph{Matching of triangulations}
The numbers of elements of $\calM_\Sigma$ not already recognized to belong
to $\calO^\hyp_\Sigma$ are as described in the central column of Table~\ref{bd:non-hyp:tab}.
As already remarked, all these manifolds come with a triangulation consisting of 4 tetrahedra.
Now, one of the features of
``Orb'' is to compare two triangulated manifolds for equality by randomizing the initial
triangulations and matching. So we have first exploited this feature to reduce the numbers of
potentially distinct homeomorphism types, getting the results described in the
right column of Table~\ref{bd:non-hyp:tab}.
In the rest of this section we describe the proof of the following result:

\begin{table}
\begin{center}
\begin{tabular}{c|c|c}
type according & apparently & apparently distinct \\
to the boundary & non-hyperbolic & after matching \\ \hline\hline
 $\calM_\emptyset$ & 37 & 17 \\ \hline
 $\calM_{\Sigma_1}$ & 70 & 21 \\ \hline
 $\calM_{\Sigma_1\sqcup \Sigma_1}$ & 7 & 5 \\ \hline
 $\calM_{\Sigma_2}$ & 50 & 16 \\ \hline
 $\calM_{\Sigma_2\sqcup \Sigma_1}$ & -- & -- \\ \hline
 $\calM_{\Sigma_3}$ & -- & -- \\
\end{tabular}
\end{center}
\mycap{Numbers of apparently non-hyperbolic
elements of $\calM$, and potentially distinct homeomorphism types
after the triangulation matching performed using ``Orb''\label{bd:non-hyp:tab}}
\end{table}

\begin{prop}\label{conclusion:prop}
For $\Sigma=\emptyset,\Sigma_1,\Sigma_1\sqcup\Sigma_1,\Sigma_2$ and
$I=17,21,5,16$, respectively, let $\big(M_\Sigma^{(i)}\big)_{i=1}^I$
be the manifolds as in the right column of
Table~\ref{bd:non-hyp:tab}. Then:
\begin{enumerate}
\item If $i\neq j$ then $M_\Sigma^{(i)}$ is not homeomorphic to
$M_\Sigma^{(j)}$;
\item Each $M_\Sigma^{(i)}$ is non-hyperbolic.
\end{enumerate}
\end{prop}

This implies Propositions~\ref{hyp:Sigma2:prop} and~\ref{hyp:Sigma1:prop},
the equalities $\calO^\non_\Sigma=\big(M_\Sigma^{(i)}\big)_{i=1}^I$ for all four relevant $\Sigma$'s,
and hence Theorem~\ref{main:teo}.
Our proof utilizes computers and theoretical work.
Note that
Proposition~\ref{conclusion:prop}
shows that ``Orb'' was totally efficient both in constructing the hyperbolic structures
and in comparing the non-hyperbolic manifolds for homeomorphism.

In the sequel we freely use several classical notions, results and
techniques of 3-manifold topology, in particular the definition of
essential surface, the Haken-Kneser-Milnor decomposition along
spheres, the definition and properties of Seifert fibred spaces, and
the Jaco-Shalen-Johansson decomposition along tori and annuli,
see~\cite{Hempel,mafo,matbook}. Moreover we use the fact that if a
manifold contains a properly embedded essential surface with
non-negative Euler characteristic then the manifold cannot be
hyperbolic.

\paragraph{The ``3-Manifold Recognizer''}
As already mentioned, besides ``Orb'' we have employed another software,
namely the ``3-Manifold Recognizer,'' written by Tarkaev and
Matveev~\cite{recognizer}. The input to this program is a
triangulation of a 3-manifold $M$
and its output is the ``name'' of $M$, by which we mean the following:
\begin{itemize}
\item For a Seifert $M$, (one of) its Seifert structure(s);
\item For a hyperbolic $M$,
its presentation(s) as a Dehn filling of a manifold in the
tables of Weeks~\cite{CaHiWe};
\item For an irreducible $M$ having JSJ decomposition into more than
one block, the names (as just illustrated) of the blocks, together
with the gluing instructions between the blocks;
\item For a reducible manifold, the names (as just illustrated) of its
irreducible summands.
\end{itemize}
The program is not guaranteed to always find the name of the
manifold (for instance, it does not even attempt to do this for
manifolds with boundary of genus 2 or more, and it happens to fail
also in other cases). But it can always compute the first homology
and, in the case of boundary of genus at most $1$, the Turaev-Viro
invariants~\cite{TV}, which turned out to be very useful for us.

We now describe the proof of Proposition~\ref{conclusion:prop},
breaking it into separate paragraphs according to the boundary type $\Sigma$, and
at the same time we provide detailed topological information on the
manifolds $M_\Sigma^{(i)}$.

\paragraph{Closed manifolds}
Let us start with the case $\Sigma=\emptyset$.
The second item in Proposition~\ref{conclusion:prop}, namely the proof
that each $M_\emptyset^{(i)}$ is non-hyperbolic, was not an issue in this case.
In fact, it has been known for a long time~\cite{matbook} that any
triangulation of a closed hyperbolic manifold contains at least 9
tetrahedra, whereas each $M_\emptyset^{(i)}$ admits a triangulation with 4 tetrahedra.

To show that $M_\emptyset^{(i)}\not\cong M_\emptyset^{(j)}$ for $1\leqslant i<j\leqslant 17$
we have run the ``Recognizer,'' that successfully identified all the manifolds
(this was also independently done by Tarkaev). From the names (all
manifolds turned out to be Seifert or connected sums of Seifert) we
could see that the $M_\emptyset^{(i)}$'s were indeed all distinct, except possibly
for $M_\emptyset^{(1)}$ and $M_\emptyset^{(2)}$, that were both recognized to be
the connected sum of two copies of the lens space
$L(3,1)$. Since $L(3,1)$ has no orientation-reversing
automorphism, even if one looks (as we do) at orientable but unoriented
manifolds, there are two distinct ways of performing the connected
sum of $L(3,1)$ with itself, so the names of
$M_\emptyset^{(1)}$ and $M_\emptyset^{(2)}$ provided by the
``Recognizer'' were indeed ambiguous.

To show that $M_\emptyset^{(1)}\not\cong M_\emptyset^{(2)}$ we then
had to examine their triangulations by hand, introducing an
arbitrary orientation on each and finding the essential sphere
realizing the connected sum. Cutting along this sphere and capping
off, we saw that for $M_\emptyset^{(1)}$ the two connected summands
were distinctly oriented copies of $L(3,1)$, while for
$M_\emptyset^{(2)}$ they were consistently oriented. This led us to
the proof of Proposition~\ref{conclusion:prop} for
$\Sigma=\emptyset$. More precisely we established the following:

\begin{prop}
The set $\mathcal{O}^\non_\emptyset$ consists of $13$ irreducible
manifolds and $4$ reducible ones. The irreducible manifolds are
the Seifert spaces
$$\matS^3,\qquad \matP^3,\qquad \matS^2\times\matS^1,\qquad L(3,1),\qquad
L(4,1),\qquad L(5,2),$$
$$L(6,1),\qquad L(9,2),\qquad L(12,5),\qquad \big(\matP^2;(3,2),(1,0)\big),$$
$$\big(\matP^2;(2,1),(1,1)\big),\qquad \big(\matP^2;(1,3)\big),
\qquad \big(\matS^2;(2,1),(3,1),(3,1),(1,-1)\big),$$
and the reducible ones are
$$\matP^3\#\matP^3,\qquad \matP^3\# L(3,1),\qquad L(3,1)\#
L(3,1),\qquad L(3,1)\#\big(-L(3,1)\big).$$
\end{prop}

\paragraph{One-cusped manifolds}
In this case both items in Proposition~\ref{conclusion:prop}
required some work. We proceeded as follows.

To prove that $M_{\Sigma_1}^{(i)}\not\cong M_{\Sigma_1}^{(j)}$
for $1\leqslant i<j\leqslant 21$
 we again employed the
``Recognizer'', using which we calculated the first homology
group and Turaev-Viro invariants up to order~16 of each
$M_{\Sigma_1}^{(i)}$. From this
computation we deduced that $M_{\Sigma_1}^{(i)}\not\cong M_{\Sigma_1}^{(j)}$
for $1\leqslant i<j\leqslant 21$ except possibly for $i=1,2,3,4$ and $j=i+4$. For
the four pairs of manifolds left, we showed the homeomorphism was
impossible by analyzing the JSJ decompositions. Specifically, $M_{\Sigma_1}^{(1)}$
and $M_{\Sigma_1}^{(5)}$ turned out to be Seifert and distinct, and the same
happened for $M_{\Sigma_1}^{(2)}$ and $M_{\Sigma_1}^{(6)}$, whereas $M_{\Sigma_1}^{(3)}$ and $M_{\Sigma_1}^{(7)}$ had
non-trivial JSJ decompositions, with the same blocks but different
gluing matrices, and analogously for $M_{\Sigma_1}^{(4)}$ and $M_{\Sigma_1}^{(8)}$.

The results just described allowed us to conclude that
$M_{\Sigma_1}^{(i)}$ is non-hyper\-bolic for $i=1,\ldots,8$. To show
that the same holds for $i=9,\ldots,21$
we used the ``Recognizer'' again to compute connected sum and JSJ decompositions.
In each instance the desired result was returned
because we obtained either connected sums or manifolds having JSJ
decomposition consisting of Seifert pieces (sometimes only one of
them). It is perhaps worth mentioning that in one case the
``Recognizer'' failed to return the answer right away, but we were
able to transform the triangulation by hand into one that the
``Recognizer'' could handle.

These arguments led us to the proof of
Proposition~\ref{conclusion:prop} for the case $\Sigma=\Sigma_1$,
and also to the next more specific result. In its statement
we use matrices to encode gluings between boundary components of
Seifert spaces, which requires choosing homology bases; when the
base surface of the fibration is orientable, the homology basis is
$(\mu,\lambda)$, where $\mu$ is a boundary component of the base
surface of the fibration and $\lambda$ is a fibre;
see~\cite{mafo} for the non-orientable case.

\begin{prop}
The $21$ elements of the set $\calO^\non_{\Sigma_1}$ subdivide as follows:
\begin{itemize}
\item $2$ reducible manifolds, both being the connected sum of two Seifert spaces;
\item $10$ irreducible Seifert spaces;
\item $7$ irreducible manifolds whose JSJ decomposition consists of two Seifert blocks;
\item $2$ irreducible manifolds whose JSJ decomposition consists of three Seifert blocks.
\end{itemize}
More precisely:
\begin{itemize}
\item The $2$ reducible manifolds are $\matP^3\#(D^2\times\matS^1)$
and  $L(3,1)\#(D^2\times\matS^1)$;
\item The $10$ Seifert spaces are
$$\begin{array}{ll}
D^2\times\matS^1,&
\big(\matS^2\setminus 3D^2,(1,0)\big),\\
\big(D^2,(2,1),(2,1),(1,0)\big),&
\big(D^2,(2,1),(3,1),(1,-1)\big),\\
\big(D^2,(3,1),(3,2),(1,0)\big),&
\big(D^2,(3,2),(3,2),(1,-1)\big),\\
\big(D^2,(3,2),(4,1),(1,-1)\big),&
\big(D^2,(3,1),(4,1),(1,0)\big), \\
\big(\matP^2\setminus D^2,(2,1),(1,0)\big),&
\big(\matP^2\setminus D^2,(3,2),(1,0)\big);
\end{array}$$
\item The $7$ manifolds having JSJ decomposition consisting of two Seifert blocks
are obtained by gluing the following pairs of Seifert spaces along the homeomorphism
represented by the matrix $\tiny\left(\begin{array}{cc}0 & 1\\
1 & 0\end{array}\right)$:
$$\begin{array}{rcl}
\big(\matS^2\setminus 2D^2,(2,1),(1,0)\big)& {\rm and} &\big(D^2,(2,1),(2,1),(1,0)\big),\\
\big(\matS^2\setminus 2D^2,(2,1),(1,1)\big)& {\rm and} &\big(D^2,(2,1),(4,3),(1,-1)\big),\\
\big(\matS^2\setminus 2D^2,(3,1),(1,-1)\big)& {\rm and} &\big(D^2,(2,1),(3,2),(1,-1)\big),\\
\big(\matS^2\setminus 2D^2,(3,2),(1,0)\big)& {\rm and} &\big(D^2,(2,1),(3,2),(1,-1)\big),\\
\big(\matS^2\setminus 2D^2,(2,1),(1,0)\big)& {\rm and} &\big(D^2,(3,1),(3,2),(1,-1)\big),\\
\big(\matS^2\setminus 2D^2,(2,1),(1,-1)\big)& {\rm and} &\big(D^2,(3,1),(3,1),(1,-1)\big),\\
\big(\matP^2\setminus 2D^2,(1,1)\big)& {\rm and} &\big(D^2,(2,1),(3,1),(1,-1)\big);
\end{array}$$
\item The $2$ manifolds having JSJ decomposition consisting of three Seifert blocks are obtained
by gluing two Seifert spaces to two different boundary components of
$\big(\matS^2\setminus 3D^2,(1,2)\big)$.  In the first example the remaining two Seifert blocks are both
$(D^2,(2,1),(3,2),(1,-1))$.  In the second example the two remaining two Seifert blocks are
$(D^2,(2,1),(3,1),(1,-1))$ and $(D^2,(2,1),(3,2),(1,-1))$.
The gluing homeomorphisms are all encoded
by the matrix $\tiny\left(\begin{array}{cc}0 & 1\\
1 & 0\end{array}\right)$.
\end{itemize}
\end{prop}

\begin{rem}
\emph{The fact that $M_{\Sigma_1}^{(3)}$ and $M_{\Sigma_1}^{(7)}$ have
JSJ decompositions with the same two blocks but different
gluing matrices, and analogously for $M_{\Sigma_1}^{(4)}$ and $M_{\Sigma_1}^{(8)}$,
can be recovered from the statement just given by changing some
parameters of the exceptional fibres. This allows one to get
identical presentations of some Seifert spaces but different gluing matrices.}
\end{rem}

\paragraph{Two-cusped manifolds}
For the case $\Sigma=\Sigma_1\sqcup\Sigma_1$ we had to deal with 5
manifolds, which we did using the ``Recognizer''. To show that they
are distinct we computed their Turaev-Viro invariants, which led to
the desired conclusion right away. To prove that they are not
hyperbolic we determined their JSJ decomposition, which always
turned out to consist of Seifert blocks, whence the conclusion. More
precisely we established the following:

\begin{prop}
All $5$ elements of $\calO^\non_{\Sigma_1\sqcup\Sigma_1}$ are irreducible. Three of them are
Seifert spaces and two have JSJ decomposition consisting of two Seifert blocks. The Seifert spaces are
$$\Sigma_1\times[0,1],\qquad\big(\matS^2\setminus 2D^2;(2,1),(1,-1)\big),\qquad \big(\matS^2\setminus 2D^2;(3,2),(1,1)\big);$$
the Seifert blocks for the two other manifolds are respectively
$\big(\matS^2\setminus 3D^2;(1,0)\big)$ and
$\big(D^2;(2,1),(3,1),(1,-1)\big)$, and
two copies of $\big(\matS^2\setminus 2D^2;(2,1),(1,-1)\big)$, while
the gluing is encoded by the matrix
$\tiny\left(\begin{array}{cc}0 & 1\\ 1 & 0 \end{array}\right)$
in both cases.
\end{prop}

\paragraph{Genus-2 boundary: distinguishing manifolds}
The case of genus-2 boundary was the hardest to settle, in particular
because it could not be dealt with using the
``Recognizer.'' We concentrate here on
the task of showing that $M_{\Sigma_2}^{(i)}\not\cong M_{\Sigma_2}^{(j)}$
for $1\leqslant i<j\leqslant 16$ (item 1 of Proposition~\ref{conclusion:prop}),
postponing the proof of non-hyperbolicity to another paragraph.
We proceeded as follows:
\begin{enumerate}
    \item We first analyzed (by hand) the Turaev-Viro
    invariants of each $M_{\Sigma_2}^{(i)}$. This allowed us to break down our set of
    16 manifolds into three groups of 4 manifolds, one group of 2,
    and two groups of 1, such that the manifolds in each group have
    the same Turaev-Viro invariants of all orders, while manifolds
    in different groups have a distinct Turaev-Viro invariant (of order 6
    or 7, as it turned out);
    \item Then we determined (by computer) the homology of the
    three-fold coverings of the manifolds in each group. This
    allowed us to conclude that $M_{\Sigma_2}^{(i)}\not\cong M_{\Sigma_2}^{(j)}$ for $1\leqslant i< j\leqslant 16$
    except possibly $i=1,2$ and $j=i+2$. Moreover, it was not
    difficult to show that $\pi_1(M_{\Sigma_2}^{(i)})=\pi_1(M_{\Sigma_2}^{(i+2)})$ for
    $i=1,2$ (and in fact previously we had also shown that $M_{\Sigma_2}^{(i)}$ and
    $M_{\Sigma_2}^{(i+2)}$ have the same Turaev-Viro invariants of all orders);
    \item To deal with the remaining two pairs $M_{\Sigma_2}^{(1)},M_{\Sigma_2}^{(3)}$ and
    $M_{\Sigma_2}^{(2)},M_{\Sigma_2}^{(4)}$, the strategy was to find their JSJ decompositions.
    Below we explain in some detail how this was done.
\end{enumerate}

The general idea was to switch from triangulations to the dual
viewpoint of special spines of 3-manifolds, and more generally to
\emph{simple} spines~\cite{matbook}. The reason why this was
beneficial in this case is that a special spine that contains a
2-component with embedded closure incident to two vertices (as our
spines turned out to do) admits a so-called inverse
L-move~\cite{matbook}, whose result is a simple spine of the same
manifold. In particular, this spine may contain an annulus (or
M\"obius strip) 2-component, and it frequently turns out that the
annulus transversal to the core of the annulus 2-component (or to
the boundary of the M\"obius strip) is essential. Moreover, if the
initial spine has a small number of vertices, one may hope that
after cutting along the annulus the spine breaks down into easily
identifiable pieces (for instance, polyhedra that collapse onto
graphs), in which case the annulus already constitutes the JSJ
splitting surface of the manifold in question. This is precisely the
strategy which worked in our case.

\begin{figure}
    \begin{center}
    \includegraphics[scale=0.47]{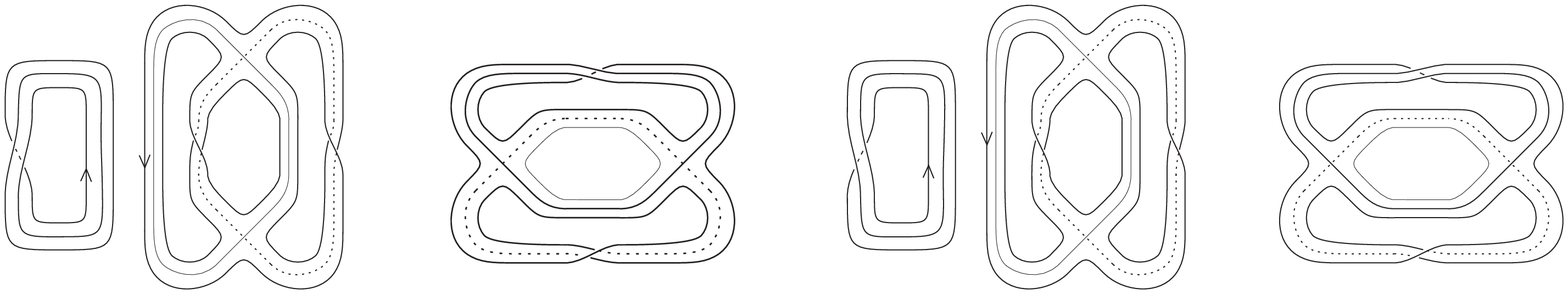}
\vspace{1cm}
    \mycap{The simple spines $P_1,\ldots,P_4$  of
$M_{\Sigma_2}^{(1)},\ldots,M_{\Sigma_2}^{(4)}$.
The picture always shows the boundary of a regular neighbourhood of the
locus of non-surface points.
To get $P_1$ from the two separate fragments shown one must
identify the two curves marked by arrows, which constitute the core $\alpha_1$
of the annular 2-component of $P_1$, while all other 2-components are discs.
The same applies to $P_3$, which contains an annulus with core $\alpha_3$.
To get $P_2$ from the fragment shown one should attach a M\"obius strip to the thick
curve $\alpha_2$ and a disc to the other one, and the same applies to $P_4$, which
contains a M\"obius strip bounded by a curve $\alpha_4$.}
    \label{spines:two:pairs:fig}
    \end{center}
    \end{figure}

Let us now turn to our specific situation. After dualizing the
triangulations and applying the inverse L-move we obtained the
simple spines $P_1,\ldots,P_4$ shown in
Fig.~\ref{spines:two:pairs:fig}.
As explained in the caption, the spines of $M_{\Sigma_2}^{(1)}$ and
$M_{\Sigma_2}^{(3)}$ contain an annular 2-component, while those of $M_{\Sigma_2}^{(2)}$ and $M_{\Sigma_2}^{(4)}$
contain a M\"obius strip 2-component. Let us denote by $S_i$ the
properly embedded annulus or M\"obius strip transversal to the
curve $\alpha_i$ also described in the
caption of Fig.~\ref{spines:two:pairs:fig}.

\bigskip

We begin with the case $i=1,3$.
As one sees from the picture, cutting
$P_i$ along $\alpha_i$ one gets a
disjoint union of two polyhedra that collapse respectively onto a circle
and onto a graph of Euler characteristic $-1$. Since this
corresponds to cutting $M_{\Sigma_2}^{(i)}$ along $S_i$, we deduce
that $M_{\Sigma_2}^{(i)}$ is obtained by gluing a genus-2 handlebody and a solid
torus along a boundary annulus. Looking at the core curves
of the glued annuli, it is not difficult to
show that the annulus $S_i$ is essential in $M_{\Sigma_2}^{(i)}$, so it gives the
JSJ decomposition. Finally, taking a closer examination of the gluings, we saw
that the annuli used in both gluings are the same, while the gluing
homeomorphisms are different. This allowed us to conclude that
$M_{\Sigma_2}^{(1)}\not\cong M_{\Sigma_2}^{(3)}$.

\bigskip

\begin{figure}
    \begin{center}
    \includegraphics[scale=0.8]{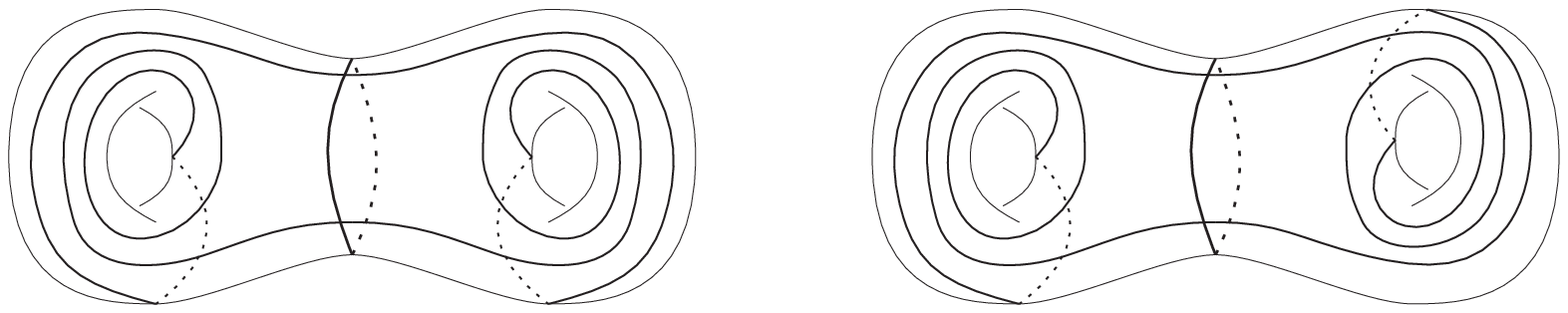}
\vspace{1cm}
    \mycap{The core curves of the annuli used to reconstruct
    $M_{\Sigma_2}^{(2)}$ and $M_{\Sigma_2}^{(4)}$.}
    \label{two:curves:fig}
    \end{center}
    \end{figure}

Let us now turn to the case $i=2,4$. Cutting $P_i$ along the core
circle of the M\"obius strip component (which again corresponds to
cutting $M_{\Sigma_2}^{(i)}$ along $S_i$) yields a polyhedron which
collapses onto a graph of Euler characteristic $-1$. Even if we get
a single polyhedron, (which must be the case since this time the cut
is along the core of a M\"obius strip), we again conclude that the
initial manifold is obtained by gluing a genus-2 handlebody and a
solid torus along a boundary annulus. As before it is not hard to
show that the annulus is in fact essential, so it gives the JSJ
decomposition. In addition, we have proved that the annulus in the
boundary of the solid torus is the same in both cases, its core
being the curve of type (2,1). On the contrary, the cores of the
annuli on the boundary of the genus-2 handlebody used to obtain
$M_{\Sigma_2}^{(2)}$ and $M_{\Sigma_2}^{(4)}$ are those shown in
Fig.~\ref{two:curves:fig}.

The conclusion that $M_{\Sigma_2}^{(2)}\not\cong M_{\Sigma_2}^{(4)}$
now follows from the next result, the long proof of which we only outline:

\begin{prop}
No homeomorphism of the genus-$2$ handlebody $H$ takes the curve $\ell_2$ shown in
Fig.~\ref{two:curves:fig}-left to the curve $\ell_4$ shown in
Fig.~\ref{two:curves:fig}-right.
\end{prop}

\begin{proof}
As already mentioned, we restrict
ourselves to indicating the general scheme of our argument only.
As one sees from Fig.~\ref{two:curves:fig}, for $i=2,4$
there exists an essential disc $D_i$ in $H$ which intersects $\ell_i$ transversely
in exactly two points. Moreover cutting $H$ along $D_i$ we get
two solid tori $T^0_i$ and $T^1_i$ such that $\partial T^j_i$ contains
a distinguished disc $\Delta^j_i$ and an arc $\beta_i^j$ properly
embedded in $\partial T^j_i\setminus\Delta^j_i$. The pair
$(H,\ell_i)$ is obtained by gluing $T^0_i$ to $T^1_i$ along a homeomorphism
$\Delta^0_i\to\Delta^1_i$, with $\ell_i$ being the image of
$\beta^0_i\cup\beta^1_i$. It is actually quite easy to see that the four triples
$(T^j_i,\Delta^j_i,\beta^j_i)$ for $i=2,4$ and $j=0,1$ can be identified to each other, but
after doing this the gluing homeomorphisms $\Delta^0_2\to\Delta^1_2$ and
$\Delta^0_4\to\Delta^1_4$ differ by a rotation of angle $\pi$,
which is isotopic to the identity but not in a way that
preserves the endpoints of the arcs. The proof of the proposition
then follows from the next:

\medskip

\noindent\textsc{Claim}. For $\ell\in\{\ell_2,\ell_4\}$, the
disc $D$ properly embedded in $H$ which intersects $\ell$ transversely in two points
and splits $H$ into two solid tori is unique up to isotopy preserving $\ell$.

\medskip

The proof of this claim is rather long and technical. We consider a
handle decomposition of $H$ into one 0-handle and two 1-handles.
This yields a decomposition of $\partial H$ into three punctured
discs, namely one sphere with four holes and two annuli. Slightly
modifying the definition in~\cite{matbook} we then call
\emph{normal} with respect to this decomposition a curve in
$\partial H$ which intersects each of the punctured discs along a
collection of simple arcs with endpoints on different boundary
components or along a simple closed curve. We next establish the
following two facts:
\begin{enumerate}
    \item Up to isotopy preserving $\ell$ there is a unique normal curve that intersects $\ell$ in two
    points and decomposes $H$ into two solid tori;
    \item The boundary of $D$ can be
    isotoped (preserving $\ell$) to normal position.
\end{enumerate}
This concludes our argument.
\end{proof}

\paragraph{Genus-2 boundary: non-hyperbolicity} To show that none of
the manifolds $M_{\Sigma_2}^{(i)}$ is hyperbolic, we used
again the idea described above. Namely, we constructed for each
$M_{\Sigma_2}^{(i)}$ a simple spine
with an annulus or M\"obius strip component and we
proved that the corresponding proper annulus in the manifold is
essential. This was done as follows:
\begin{enumerate}
    \item For about half of the $M_{\Sigma_2}^{(i)}$'s, the special spine
    dual to the initial triangulation already contained
    a 2-component incident to two vertices, so we found a simple
    spine with an annulus or M\"obius strip 2-component by
    applying an inverse L-move, as above. For the other $M_{\Sigma_2}^{(i)}$'s
    we did the same but we first had to change the initial special spine, by
    applying first one positive
    $T$-move~\cite{matbook} and then one inverse $T$-move elsewhere.
    \item From the spine of $M_{\Sigma_2}^{(i)}$ constructed
    in the previous item we got a properly embedded annulus $S_i$, that we then
    showed to be essential.
We did this by cutting
$M_{\Sigma_2}^{(i)}$
along $S_i$, which gave the following:
\begin{enumerate}
\item In $2$ cases, a genus-2 handlebody;
\item In $6$ cases,
the union of a genus-2 handlebody and a solid torus;
\item In $4$ cases, a manifold that could be further split along an annulus
into the union of a genus-2 handlebody and a solid torus;
\item In $4$ cases, the union of a solid torus and a manifold
as described in the previous point.
\end{enumerate}
In all cases, analyzing the way $M_{\Sigma_2}^{(i)}$ can be reconstructed from the
pieces $S_i$ cuts into, we could then show that it is irreducible and that
within it
$S_i$ is $\pi_1$-injective and
not boundary-parallel, from which we got the desired conclusion.
\end{enumerate}

\paragraph{Further information for genus-2 boundary} The decomposition
(a)-(d) just described along annuli of the 16 elements of
$\calO^\non_{\Sigma_2}$ provides a rather accurate description of the
topology of these manifolds. In addition to it, we mention that in cases
(c) and (d) the second splitting annulus is not disjoint from the trace of
$S_i$, so the splitting cannot be described as being along the union of
two disjoint annuli.

\vspace{1cm}

\noindent RedTribe\\
Carlton, Victoria\\
Australia 3053\\
{\tt damianh@redtribe.com}

\vspace{.5cm}

\noindent Dipartimento di Matematica Applicata\\
Via Filippo Buonarroti, 1C\\
56127 PISA -- Italy\\
{\tt petronio@dm.unipi.it}

\vspace{.5cm}

\noindent Dipartimento di Matematica Applicata\\
Via Filippo Buonarroti, 1C\\
56127 PISA -- Italy\\
{\tt pervova@csu.ru}

\end{document}